\documentclass[11pt]{amsart}
\usepackage{graphicx}

\usepackage{amssymb}
\usepackage{amsthm}
\usepackage{latexsym}
\newcommand{\sub}{\subseteq}

\newcommand{\lra}{\Leftrightarrow}

\newcommand{\ra}{\Rightarrow}

\newcommand{\sm}{\setminus}

\newcommand{\al}{\alpha}
\newcommand{\La}{\Lambda}
\newcommand{\be}{\beta}

\newcommand{\mc}{\mathcal}

\newcommand{\F}{\mc F}

\newcommand{\op}{\operatorname}
\newcommand{\tmax}{t\op{-Max}}
\newcommand{\tspec}{t\op{-Spec}}
\newcommand{\spec}{\op{Spec}}

\newcommand{\fbar}{\boldsymbol{\overline{F}}}

\newtheorem{thm}{Theorem}[section]
\newtheorem{prop}[thm]{Proposition}
\newtheorem{lemma}[thm]{Lemma}
\newtheorem{cor}[thm]{Corollary}
\theoremstyle{definition}

\newtheorem{rem}[thm]{Remark}

\begin{document}

\title{Star stable domains}%
\author{Stefania Gabelli}
\address{Dipartimento di Matematica, Universit\`a
degli Studi ``Roma Tre'',
Largo San Leonardo Murialdo 1, 00146 Roma}%
\email{gabelli@matrm3.mat.uniroma3.it, picozza@matrm3.mat.uniroma3.it}%

\author{Giampaolo Picozza}

\email{}

\date{\today}




\begin{abstract}
We introduce and study the notion of $\star$-stability with respect to a semistar operation $\star$ defined on a domain $R$; in particular we consider the case where $\star$ is the $w$-operation.
This notion allows us to generalize and improve several properties of stable domains and totally divisorial domains.
\end{abstract}

\maketitle

\section*{Introduction}

Star operations, as the $v$-closure (or divisorial closure), the $t$-closure and the $w$-closure are
an essential tool in modern multiplicative ideal theory for characterizing and investigating several classes of integral domains.
For example, in the last few decades a large amount of literature has appeared on  \emph{Mori domains}, that is domains satisfying the ascending chain condition on divisorial ideals, and \emph{Pr\"ufer $v$-multiplication domains}, for short \emph{P$v$MDs},  that is domains in which each finitely generated ideal is $t$-invertible (or $w$-invertible). The consideration that some important operations on ideals, like the integral closure, satisfy almost all the properties of star operations led A. Okabe and R. Matsuda to introduce in 1994 the more general and flexible notion of semistar operation \cite{OM94}. The class of semistar operations includes the classical star operations and often provides a more appropriate context for approaching several questions of multiplicative ideal theory, see for example \cite{EFP, FH2000, fjs, fl01, fl03, P}.
In this paper, we introduce the notion of $\star$-stability with respect to a semistar operation $\star$.

Motivated by earlier work of H. Bass \cite{b} an J. Lipman
\cite{L} on the number of generators of an ideal, in 1974
 J. Sally and W. Vasconcelos defined a Noetherian ring $R$ to be
\emph{stable} if each nonzero ideal of $R$ is projective over its
endomorphism ring End$_R(I)$ \cite{sv}. In a note of 1987, D.D. Anderson, J. Huckaba and
I. Papick  considered the notion of stability for
arbitrary integral domains \cite{AHP}. When $I$ is a nonzero ideal of a domain $R$,
then End$_R(I)=(I:I)$; thus a domain $R$ is stable if each nonzero
ideal $I$ of $R$ is invertible in the overring $(I:I)$. Since 1998,  stable
domains have been thoroughly investigated by
B. Olberding in a series of papers \cite{olb98, olb, olb01, olb01-2, olb02}.

Given a semistar operation $\star$ on a domain
$R$, we say that a nonzero ideal $I$ of $R$ is \emph{$\star$-stable} if
$I^\star$ is $\dot{\star}$-invertible in $(I^\star:I^\star)$ and
that $R$ is \emph{$\star$-stable} if each nonzero ideal of $R$ is
$\star$-stable. (Here we denote by $\dot{\star}$ the
semistar operation induced by $\star$ on a fixed overring $T$ of $R$.) This notion allows
us to generalize and improve several properties of stable domains and totally divisorial domains.  We also recover some results proven in \cite[Section 2]{SS2} for $\star=w$.

Even though many results are stated for a general semistar operation, for technical reasons, the most interesting consequences are obtained for (semi)star operations spectral and of finite type. In this case, we show that $\star$-stability implies that $\star$ is the $w$-operation on $R$; in particular, on stable domains the $w$-operation is the identity.

For a (semi)star operation spectral and of finite type, the main result of Section 1 is that a domain $R$ is $\star$-stable if and only if $R$
is $\star$-locally stable and has $\star$-finite character, if and only if $R$
is $\star$-locally stable and each $\star$-ideal of $R$ is $\dot \star$-finite in its endomorphism ring.
This implies that if a domain is locally stable, then stability is equivalent to the property that each nonzero ideal $I$ is finitely generated in the overring $(I:I)$.

In Section 2 we study $\star$-stability of overrings and we show that, for semistar operations of finite type,  the $\star$-integral closure of a $\star$-stable domain is a P$v$MD.

In Section 3 we extend some properties of totally divisorial domains in the setting of semistar operations. For $\star=w$, we prove that each $t$-linked overring $T$ of $R$ is $\dot{w}$-divisorial if and only if  all the endomorphism rings of $w$-ideals are $\dot{w}$-divisorial, if and only if $R$ is $w$-stable and $w$-divisorial. Under these conditions, $\dot{w}$ is the $w$-operation on $T$.
As a consequence, we get that  $R$ is totally divisorial if and only if all the overrings of type $(I:I)$ are divisorial, if and only if each nonzero ideal $I$ of $R$ is $m$-canonical in $(I:I)$.
The Mori case and the integrally closed case are of particular interest.

Finally, in Section 4 we show that $w$-stable $w$-divisorial domains are $v$-coherent and use this fact to show that $w$-stable $w$-divisorial (respectively, totally divisorial) domains  share several properties with generalized Krull (respectively, Dedekind) domains. As a matter of fact, in the integrally closed case each one of these properties becomes equivalent  to $R$ being a generalized Krull (respectively, Dedekind) domain; so that a $w$-stable $w$-divisorial (respectively, totally divisorial) domain can be viewed as a \textquoteleft\textquoteleft non-integrally closed generalized Krull (respectively, Dedekind) domain".

\medskip
Throughout this paper $R$ will be an integral domain with quotient field $K$, $R\neq K$. 
We denote by $\boldsymbol{F}(R)$ the set of nonzero
fractional ideals of $R$, by $\fbar(R)$ the set of nonzero
$R$-submodules of $K$ and by $\boldsymbol{f}(R)$ the set of
nonzero finitely generated
$R$-submodule of $K$. Clearly $\boldsymbol{f}(R) \sub \boldsymbol{F}(R) \sub \fbar(R)$.

A \emph{semistar operation} on $R$ is a map $\star:
\fbar(R) \to \fbar(R)$ such that, for each $E,F \in \fbar(R)$ and
for each $x \in K$, $x \neq 0$, the following properties hold:

\begin{enumerate}
\item[$(\star_1)$] $(xE)^\star=xE^\star$;
\item[$(\star_2)$] $E
\subseteq F$ implies $E^\star \subseteq F^\star$;
\item[$(\star_3)$] $E \subseteq E^\star$ and $E^{\star \star}:=
\left(E^\star \right)^\star=E^\star$.
\end{enumerate}

Recall that, for all $E,F \in \boldsymbol{\overline{F}}(R)$, we
have:
$$\begin{array}{rl} (EF)^\star =& \hskip -7pt (E^\star
F)^\star =\left(EF^\star\right)^\star
=\left(E^\star F^\star\right)^\star\,;\\
(E+F)^\star =& \hskip -7pt \left(E^\star + F\right)^\star= \left(E
+
F^\star\right)^\star= \left(E^\star + F^\star\right)^\star\,;\\
(E:F)^\star \subseteq & \hskip -7pt (E^\star:F^\star) = (E^\star
:F) =
\left(E^\star:F\right)^\star\,;\\
(E\cap F)^\star \subseteq & \hskip -7pt E^\star \cap F^\star =
\left(E^\star \cap F^\star \right)^\star,\,\mbox{ if $E\cap F
\neq (0)$}\,;
\end{array}
$$
see for instance \cite{FH2000}.

If $\star_1$ and $\star_2$ are semistar operations on $R$, we say
that $\star_1 \leq \star_2$ if $E^{\star_1} \subseteq
E^{\star_2}$, for each $E \in \fbar(R)$. This is equivalent to the
condition that $(E^{\star_1})^{\star_2}= (E^{\star_2})^{\star_1} =
E^{\star_2}$, for each $E \in \fbar(R)$.

The identity is a semistar operation, denoted by $d$. It
follows from $(\star_3)$ that $d\leq \star$, for each semistar
operation $\star$.

A semistar operation $\star$ is
called a \emph{semistar operation of finite type} if, for each $E
\in \fbar(R)$, we have
$$E^\star = \bigcup \{F^\star \, \vert \, F
\in \boldsymbol{f}(R) \mbox{ and } F \subseteq E\}.$$

If $\star$ is any semistar operation, the semistar operation $\star_f$ defined by
$$E^{\star_{f}}:=\bigcup \{F^\star \, \vert \, F \in
\boldsymbol{f}(R) \mbox{ and } F \subseteq E\},$$
for each $E \in \boldsymbol{\overline{F}}(R)$, is
a semistar operation of finite type and $\star_f \leq \star$.

A nonzero ideal $I$ of $R$ is $\star$-finite if there exists a
finitely generated $J$ such that $ I^\star=J^\star =J^{\star_f}$.

When $R^\star=R$, $\star$ is called a \emph{(semi)star operation}
on $R$ and its restriction to the set of nonzero fractional ideals $\boldsymbol{F}(R)$
is a star operation, still denoted by $\star$.

As usual, we denote by $v$ the (semi)star operation defined by
$E^v:= (R:(R:E))$, for each $E\in \fbar(R)$, and set $t:= v_f$.
As a star operation on $R$, $v$ is called the \emph{divisorial
closure}. It is well known that  $\star \leq v$ and
$\star_f \leq t$, for each (semi)star operation $\star$
\cite[Proposition 1.6]{FH2000}.

We say that a nonzero ideal $I$ of $R$ is a
\emph{quasi-$\star$-ideal} if $I^\star \cap R =I$. A
\emph{quasi-$\star$-prime (ideal)} is a prime quasi-$\star$-ideal
and a \emph{quasi-$\star$-maximal ideal} is a quasi-$\star$-ideal
maximal in the set of all proper quasi-$\star$-ideals. A
quasi-$\star$-maximal ideal is a prime ideal \cite[Lemma
4.20]{FH2000} and, when $\star$ is a semistar operation of finite
type, each quasi-$\star$-ideal is contained in a
quasi-$\star$-maximal ideal \cite[Lemma 4.20]{FH2000}. The set of
quasi-$\star_f$-maximal ideals of $R$ will be denoted by $\star_f
\op{-Max}(R)$.
  We say that $R$ has \emph{$\star_f$-finite character}
if each nonzero ideal of $R$ is contained at most in a finite
number of quasi-$\star_f$-maximal ideals of $R$.

When $\star$ is a (semi)star operation, an ideal $I$ is a
quasi-$\star$-ideal if and only if $I^\star=I$. In this case,
like in the classical case of star operations,
we say that $I$ is a
\emph{$\star$-ideal} and, analogously, we call a
quasi-$\star$-prime ideal a \emph{$\star$-prime} and a
quasi-$\star$-maximal ideal
a \emph{$\star$-maximal ideal}.
A $v$-ideal of $R$ is also called a \emph{divisorial ideal}.

If $\star$ is a semistar operation on $R$, we denote by
$\tilde{\star}$ the semistar operation defined by
\[E^{\tilde{\star}}:= \bigcap_{M\in \star_f \op{-Max}(R)}ER_M =
\bigcup_{F\in \boldsymbol{f}(R),\,F^{\star_f}=R} (E:F),\]
for each $E\in\fbar(R)$.
We have $I^{\tilde{\star}}R_M = IR_M$, for each nonzero ideal $I$ of $R$
and each quasi-$\star_f$-maximal ideal $M$ of $R$ \cite[Lemma
4.1(2)]{FH2000}.
Clearly $\tilde{\star}=\widetilde{\star_f}$.

The semistar operation $\tilde{\star}$ is of finite type and
spectral (a semistar operation $\star$ is \emph{spectral} if there exists
$\La \subseteq \spec(R)$ such that $E^\star = \bigcap \{ER_P\,\vert\, P\in
\La\}$, for each $E \in \fbar(R)$).
More precisely, $\star
= \tilde{\star}$ if and only if $\star$ is spectral and of finite
type, if and only if $\star$ is of finite type and $(E\cap
F)^\star = E^ \star \cap F^\star$, for $E, F \in \fbar(R)$ such that $E\cap F
\neq (0)$
\cite[Corollary 3.9 and Proposition 4.23]{FH2000}.

Always we have $\tilde{\star} \leq \star_f \leq \star$.
In addition, setting $w:= \tilde{v}$, if $\star$ is a (semi)star operation,
we have $\tilde{\star} \leq w$.

A nonzero ideal $I$ of $R$ is \emph{$\star$-invertible} if
$(I(R:I))^\star = R^\star$. When $\star = \star_f$, this is
equivalent to the fact that $I(R:I)$ is not contained in any
quasi-$\star_f$-maximal ideal. Since quasi-$\star_f$-maximal
ideals and quasi-$\tilde{\star}$-maximal ideals coincide
\cite[Corollary 3.5(2)]{fl03}, it follows that an ideal $I$ is
$\star_f$-invertible if and only if it is
$\tilde{\star}$-invertible. When $\star = \star_f$, a
$\star$-invertible ideal is $\star$-finite.

If $\star$ is a semistar operation on $R$ and $T$ is an overring of $R$,
 the restriction of $\star$ to the set of $T$-submodules of $K$ is a semistar operation on $T$,  here denoted by $\dot{\star}$. When $T^ \star =T$,  $\dot{\star}$ is a (semi)star operation on $T$
\cite[Proposition 2.8]{fl01}.
Note that $\dot{\star}$ shares many properties with $\star$ (see
for instance \cite[Proposition 3.1]{giampa}); for example,
if $\star$ is of finite type then $\dot{\star}$ is of
finite type \cite[Proposition 2.8]{fl01}.

\section{$\star$-stable domains}

Let $R$ be an integral domain and $\star$ a semistar operation on
$R$. Given a nonzero fractional ideal $I$ of $R$, consider the overring
$T:=(I^\star:I^\star)$ of $R$. It is easy to see that,
$T=T^\star$; hence the restriction
of $\star$ to the set of the $T$-submodules of $K$
is a (semi)star operation on $T$, denoted by
$\dot{\star}$.

We say that
a nonzero fractional ideal $I$ of $R$ is \emph{$\star$-stable} if $I^\star$ is
$\dot{\star}$-invertible in $(I^\star:I^\star)$ and that $R$ itself is
\emph{$\star$-stable} if each nonzero (fractional) ideal of $R$ is $\star$-stable.
The notion of $d$-stable domain
coincides with the notion of stable domain introduced
in \cite{sv}.

A $\star$-invertible ideal $I$ of $R$ is $\star$-stable. In fact, since $R^\star\sub T^\star=T$,
if $(I(R:I))^\star=R^\star$, we have
$$T= R^\star T \sub ((I(R:I))^\star T)^\star= (I(R:I)T)^\star\sub(I^\star(T:I^\star))^\star\sub T$$
Thus $(I^\star(T:I^\star))^\star = T$. It follows that a domain
with the property that each nonzero ideal is $\star$-invertible is
$\star$-stable. For example any completely integrally closed
domain, is $v$-stable. Recalling that if $I$ is $v$-invertible, we
have $R=(I^t:I^t)=(I^v:I^v)$, we see that a completely integrally
closed domain that is not a P$v$MD is  $v$-stable but not
$t$-stable. Since on Krull domains $t=v$ and any integrally closed
stable domain is a Pr\"ufer domain \cite[Proposition 2.1]{rush},
any Krull domain that is not Dedekind is $t$-stable but not
stable.

\begin{prop} \label{rpunto} Let $\star$ be a semistar operation on an
integral domain $R$ satisfying one of the following conditions:

\begin{enumerate}
\item[(1)] $(E \cap F)^\star = E^\star \cap F^\star$, for each
$E,F \in \fbar(R)$  such that $E\cap F
\neq (0)$.

\item[(2)] $(R: R^\star) \neq (0)$.
\end{enumerate}

Then $R$ is $\star$-stable if and only if $R^\star$ is
$\dot{\star}$-stable.
\end{prop}

\begin{proof}
If $I$ is a nonzero ideal of $R$, then $I^\star$ is an ideal of
$R^\star$. Hence, if $R^\star$ is $\dot{\star}$-stable,
 $R$ is $\star$-stable, without any condition on $\star$.

Conversely, let $R$ be $\star$-stable and let $J$ be a nonzero ideal of
$R^\star$. Assume that condition (1) holds and consider the ideal $I:= J \cap R$
of $R$. Then $I^\star = (J \cap R)^\star = J^\star \cap R^\star = J^\star$.
It follows that $J$ is $\star$-stable.
If (2) holds, then $J$ is a fractional ideal of $R$. Hence it is $\star$-stable.
\end{proof}

Since we will be mostly interested in the case where $\star$ is a
semistar operation spectral and of finite type, that is where $\star
=\tilde{\star}$, by the previous proposition often we will restrict
 ourselves to assume that $R=R^\star$, that is to consider
(semi)star operations.

Our first result is a generalization of \cite[Theorem 3.5(1)$\Leftrightarrow$(2)]{olb02}:

\begin{prop} \label{prop:1} The following conditions are equivalent for
an integral domain $R$ and a semistar operation $\star$ on $R$:

\begin{enumerate}
\item[(i)] $R$ is $\star$-stable.

\item[(ii)] For each nonzero ideal $I$ of $R$, $I^\star$ is a divisorial
ideal of $(I^\star:I^\star)$ (that is, $ (I^{\star})^{v^\prime} =
I^\star$, where $v^\prime$ is the $v$-operation on
$(I^\star:I^\star)$).

\end{enumerate}
\end{prop}

\begin{proof}
(i)$\Rightarrow$(ii) Since $I^\star$ is $\dot{\star}$-invertible
in $(I^\star: I^\star)$ we have $I^\star =
(I^\star)^{\dot{\star}} = (I^\star)^{v^\prime}$ \cite[Lemma
2.1]{chang/park}.

(ii)$\Rightarrow$(i) Let $T:= (I^\star: I^\star)$ and let
$J:=(T:I^\star)$. We have to show that
$(I^\star J)^{\star}=T^\star = T$.

First, we show that $(J:J) = T$. We have $(J:J) =
((T:I^\star):(T:I^\star)) = ((T:(T:I^\star)):I^\star) =
((I^\star)^{v^\prime}: I^\star) = (I^\star:I^\star)=T.$

Next step is to show that $(T:I^\star J)=T$. We have
$(T:I^\star J)
= ((T:J):I^\star) = ((T:(T:I^\star)):I^\star)=((I^\star)^{v^\prime}
: I^\star) = (I^\star:I^\star)=T.$

Now we prove that $((I^\star J)^\star:I^\star J) = T$. It is clear
that $((I^\star J)^\star:I^\star J) \supseteq T$. Conversely, if
$x \in ((I^\star J)^\star:I^\star J)$ then $x(I^\star J) \subseteq
(I^\star J)^\star \subseteq T$. Hence, $((I^\star J)^\star:I^\star
J) \subseteq (T:I^\star J) = T$.

Finally, since $((I^\star J)^\star:(I^\star J)^\star)= T$, by hypothesis,
$(I^\star J)^\star = ((I^\star J)^\star)^{v^\prime}$. Thus
$T =
(T:T) = (T: (T:I^\star J)) = (I^\star J)^{v^\prime} =
(I^\star J)^{\star},$
and $I$ is $\dot{\star}$-invertible.
\end{proof}

If $I$ is a nonzero ideal of $R$, we denote by $v(I)$ the semistar
operation defined on $R$ by $E \mapsto (I:(I:E))$, for each $E\in
\fbar(R)$. When $(I:I) = R$, then $v(I)$ is a (semi)star
operation. The ideal $I$ is called \emph{$m$-canonical} if
$J^{v(I)}:=(I:(I:J))=J$, for each nonzero fractional ideal $J$ of
$R$ \cite{mcanI}.

\begin{lemma} \label{lemma:0}
Let $R$ be an integral domain and let $I$ be a nonzero ideal of $R$.
Then:
\begin{enumerate}
\item[(1)] $I^{v(I)}=I$.

\item[(2)] If $\star$ is a semistar operation on $R$ such that
$I^\star = I$, then $\star \leq v(I)$.
\end{enumerate}
\end{lemma}
\begin{proof}
(1) This is an easy consequence of the fact that $I$ is an ideal of
$(I:I)$; thus $(I:(I:I)) = I$.\\
(2) Let $E\in \fbar(R)$. Since $I=I^\star$, then
$(I:E^\star)=(I:E)$. Hence $(E^\star)^{v(I)} = (I:(I:E^\star)) =
(I:(I:E)) = E^{v(I)}$ and so $\star \leq v(I)$.
\end{proof}

\begin{prop} \label{lemma:1} The following conditions are equivalent for an integral
domain $R$ and a (semi)star operation $\star$ on $R$:
\begin{enumerate}

\item[(i)] $R$ is $\star$-stable.

\item[(ii)] $v(I^\star) = v^\prime$, for each nonzero ideal $I$ of $R$
(where $v(I^\star)$ is defined on $(I^\star:I^\star)$ and
$v^\prime$ is the $v$-operation of $(I^\star:I^\star)$).

\item[(iii)] If $I,J$ are two nonzero ideals of $R$ such that
$(I^\star:I^\star) = (J^\star:J^\star)$ then $v(I^\star) =
v(J^\star)$ (as (semi)star operations on $(I^\star:I^\star) =
(J^\star:J^\star)$).
\end{enumerate}
\end{prop}
\begin{proof}
(i)$\Rightarrow$(ii) Since $v(I^\star)$ is a (semi)star operation on
$(I^\star:I^\star)$ we have $v(I^\star) \leq v^\prime$.
Conversely, since by Proposition
\ref{prop:1} $I^\star$ is divisorial in $(I^\star:I^\star)$, as a
consequence of Lemma \ref{lemma:0}(2), we have that $v^\prime \leq v(I^\star)$.\\
(ii)$\Rightarrow$(i) We have $(I^\star)^{v^\prime} =
(I^{\star})^{v(I^\star)} = I^\star$ by Lemma \ref{lemma:0}(1). So,
$I^\star$ is divisorial in $(I^\star:I^\star)$ for each ideal $I$
of $R$ and $R$ is $\star$-stable by Proposition \ref{prop:1}.\\
(ii)$\Rightarrow$(iii) It is straightforward since both
$v(I^\star)$ and $v(J^\star)$ coincide with $v^\prime$ in $(I^\star:I^\star) = (J^\star:J^\star)$.\\
(iii)$\Rightarrow$(ii) Note that $T:=(I^\star:I^\star)$ is a
fractional ideal of $R$, since $I^\star$ is an ideal of $R$. So,
there exists a nonzero integral ideal $J$ of $R$ and a nonzero
element $x \in R$ such that $T=x^{-1}J^\star$. Clearly, $T = (T:T)=
(J^\star:J^\star)$. Thus, by hypothesis, $v(I^\star)=v(J^\star)$.
Moreover, it is easy to see that $v(J^\star) = v(x^{-1}J^\star) =
v(T)=v^\prime$, the $v$-operation of $T$. It follows that
$v(I^\star)=v^\prime$.
\end{proof}

\begin{prop} \label{max}
Let $R$ be an integral domain and $\star$ a (semi)star operation of
finite type on $R$. If $R$ is $\star$-stable, then each
$\star$-maximal ideal of $R$ is divisorial. In particular,
$\star\op{-Max}(R) = \tmax(R) = v\op{-Max}(R)$.
\end{prop}
\begin{proof}
Let $M$ be a $\star$-maximal ideal of $R$ and suppose that $M$ is not
divisorial. Then, $M^v = R$, otherwise $M^v$ would be a
$\star$-ideal containing $M$. Hence $(M:M) = (R:M)= R$. It
follows that $M$ is $\star$-invertible in $R$ and $M^\star=(M(M:M))^\star=(M(R:M))^\star=R$, a contradiction. The
second statement follows easily.
\end{proof}

\begin{cor} \label{til=w}
Let $R$ be an integral domain and $\star$ a (semi)star operation
of finite type on $R$. If $R$ is $\star$-stable then
$\tilde{\star}=w$.
\end{cor}

\begin{rem}

(1) It is possible to prove that, given two semistar operations
$\star_{1} \leq \star_2$ on $R$ either of finite type or with the
property that $(I\cap J)^{\star_i} = I^ {\star_i} \cap
J^{\star_i}$, for any pair of nonzero ideals $I, J$ of $R$ and $i=1, 2$
(for example, two spectral semistar operations), then
$\star_1$-stability implies $\star_2$-stability. Thus, for
example, a $w$-stable domain is $t$-stable. We have no examples
of $t$-stable domains that are not $w$-stable.

(2) In general it is not true  that if $\star$ is a (semi)star operation of finite type and $R$
is $\star$-stable then $\star = t$. However, we will show in Corollary \ref{staric}
that this happens when $R$ is $\star$-integrally closed.

For an example, in \cite[Example 5.4]{sv} it is proved
that the $1$-dimensional local domain $A$ with
maximal ideal $3$-generated constructed in \cite{FR} is stable. It
is clear that $A$ is Noetherian. Hence $d\neq t$ on $A$,
 because the maximal ideal is not $2$-generated  \cite[Lemma 3.5]{olb01}.
\end{rem}

\begin{lemma} \label{lemma:2}
Let $R$ be an integral domain and $\star$ a
semistar operation on $R$. Let $J$ be a nonzero ideal of $R$  and  assume that $J^\star$ is
$\dot{\star}$-finite in $(J^\star: J^\star)$. Then, for each
prime ideal $P$ of $R$:

\begin{enumerate}
\item[(1)] $(J^\star:J^\star)R_P = (J^\star R_P:J^\star R_P)$.

\item[(2)] $((J^\star:J^\star ):J^\star)R_P = ((J^\star:J^\star)R_P:
J^\star R_P)$.

\end{enumerate}

\end{lemma}
\begin{proof}
(1) Let $T:=(J^\star:J^\star)$. Since
$J^\star$ is $\dot{\star}$-finite, there exist $x_1, x_2, \ldots,
x_n \in J^\star $, such that $J^\star = (x_1T + x_2T + \ldots +
x_nT)^{\star}$. Let $H:= x_1 R+ x_2 R + \ldots + x_nR
\subseteq J^\star$, so that $(HT)^{\star} = J^\star$. Then,
$HR_P \subseteq J^\star R_P$ and
$TR_P = (J^\star:J^\star)R_P \subseteq (J^\star
R_P: J^\star R_P) \subseteq (J^\star R_P: HR_P) =
(J^\star:(HT)^\star)R_P = (J^\star: J^\star)R_P = TR_P$.
Hence, $ (J^\star R_P:J^\star R_P) = TR_P =
(J^\star: J^\star)R_P$.

(2) Since $T^\star = T$, we have $(T:J^\star)R_P \subseteq
(TR_P:J^\star R_P) \subseteq (TR_P: HR_P) = (T:HT)R_P =
(T:(HT)^\star)R_P = (T:J^\star)R_P$. Hence, $(T:J^\star)R_P =
(TR_P: J^\star R_P)$.
\end{proof}

The next result  shows in particular that the study of $\star$-stable domains can be
reduced to the local case.

\begin{thm} \label{thm:1}
Let $R$ be an integral domain and $\star$ a (semi)star operation on
$R$. If $\star=\tilde{\star}$, the following conditions are equivalent:
\begin{enumerate}
\item[(i)] $R$ is $\star$-stable.

\item[(ii)] $R$ has $\star$-finite character and $R_M$ is stable,
for each $M\in \star\op{-Max}(R)$.

\item[(iii)] $J^\star$ is $\dot{\star}$-finite
in $(J^\star: J^\star)$,
for each nonzero ideal $J$ of $R$, and $R_M$ is stable, for each $M\in \star\op{-Max}(R)$.
\end{enumerate}

Under these conditions, $\star = w$.

\end{thm}
\begin{proof}
(i)$\Rightarrow$(ii) First, we show that $R_M$ is stable, for each $M\in \star\op{-Max}(R)$.
 Let $I$
be a nonzero ideal of $R_M$. There exists an ideal $J$ of $R$ such that
$I = JR_M = J^{\star}R_M$. Since $R$ is $\star$-stable, $J^\star$
is $\dot{\star}$-invertible in $T:=(J^\star:J^\star)$, that is,
$(J^\star(T:J^\star))^\star = (J(T:J))^\star = T$ (as usual, we
denote by $\dot{\star}$ the restriction of $\star$ to the set of
the fractional ideals of $T$). In particular, $J^\star$ is $\dot{\star}$-finite
in $T$. Hence, by Lemma
\ref{lemma:2}(2), we have $(I:I)= (JR_M:JR_M) = TR_M =
(J(T:J))^\star R_M = JR_M (T:J)R_M= JR_M (TR_M:JR_M)=I((I:I):I)$.  It follows 
that $I$ is invertible in $(I:I)$ and so $R_M$ is a stable domain.

To prove that $R$ has $\star$-finite character, we prove that a family
of $\star$-maximal ideals that has nonempty intersection is a finite
family.

Let $M$ be a $\star$-maximal ideal.
Since $R_M$ is stable, by \cite[Lemma 3.1]{olb02}, $MR_M$ is
principal in $(MR_M:MR_M)$, that is, there exists $x \in MR_M$
such that $MR_M = x(MR_M: MR_M)$ and so $M^2R_M \subseteq xMR_M
\subseteq xR_M \subseteq MR_M$. The ideal $I:= xR_M \cap R$ is a
$t$-ideal (it is the contraction of a $t$-ideal of $R_M$), and so
a $\star$-ideal, since $\star \leq t$. Moreover, $IR_M = xR_M$.
Note that $M^2 \subseteq I$ and so $IR_N = R_N$ for each
$\star$-maximal ideal $N \neq M$. Since by Lemma \ref{lemma:2}(1)
$(I:I)R_N = (IR_N: IR_N)$ for each $\star$-maximal ideal $N$, we
have $(I:I) = \bigcap \{(I:I)R_N \vert N \in \star\op{-Max}(R) \}
= \bigcap \{(IR_N:IR_N) \vert N \in \star\op{-Max}(R) \} =
\bigcap \{(R_N: R_N) \vert N \in \star\op{-Max}(R), N \neq M \}
\cap (xR_M: xR_M) = \bigcap \{(R_N: R_N) \vert N \in
\star\op{-Max}(R) \} = \bigcap \{R_N \vert N \in
\star\op{-Max}(R)\} = R$. It follows that, since $R$ is
$\star$-stable, $I$ is $\star$-invertible in $R$.

Now, let $\{M_\alpha\}$ be a collection of $\star$-maximal ideals
such that $\bigcap_\alpha M_\alpha \neq (0)$. For each $M_\alpha$
we have a $\star$-ideal $I_\alpha$, constructed as above. If $y
\in \bigcap_\alpha M_\alpha$, $y \neq 0$, then $y^2 \in
\bigcap_\alpha {M_\alpha}^2 \subseteq \bigcap_\alpha I_\alpha$.
Then, $I:= \bigcap_\alpha I_\alpha \neq (0)$. Let $J = \sum_\alpha
(R:I_\alpha)$. Since for each $\alpha$, $I_\alpha$ is a
$\star$-invertible $\star$-ideal, we have ${I_\alpha} = I_\alpha^v$
\cite[Lemma 2.1]{chang/park}, and so $(R:J) = \bigcap_\alpha
(R:(R:I_\alpha)) = \bigcap_\alpha I_\alpha = I$. Note that
$I_\alpha \not \subseteq M_\beta$ if $\beta \neq \alpha$ and that,
by Lemma \ref{lemma:2}(2), $(R:I_\alpha)R_M = (R_M:I_\alpha R_M)$.
Then, for each $\alpha$, we have $JR_{M_\alpha} = \sum_\beta
(R:I_\beta)R_{M_\alpha} = (R_{M_\alpha}:IR_{M_\alpha}) +
\sum_{\beta \neq \alpha} (R_{M_\alpha}: R_{M_\alpha}) =
(R_{M_\alpha}:IR_{M_\alpha})$. Similarly, for a $ \star$-maximal
$N \not \in \{M_\alpha\}$, we have $JR_N = R_N$. Hence
$(J^\star: J^\star) \subseteq \bigcap (JR_M:JR_M) =
\bigcap(R_{M_\alpha}: I_{\alpha}R_{M_\alpha}(R_{M_\alpha}:
I_\alpha R_{M_\alpha})) \cap (\bigcap \{R_N \vert N \in  \star\op{-Max}(R),
N \not \in \{M_\alpha\} \} ) = \bigcap R_M = R$. It follows that
$J$ is $\star$-invertible in $R$ and so $\star$-finite (and so,
$t$-finite). Then, there exists $I_{\alpha_1}, I_{\alpha_2},
\ldots, I_{\alpha_n}$ such that $J^\star = ((R:I_{\alpha_1}) +
(R:I_{\alpha_2}) + \ldots + (R:I_{\alpha_n}))^v$. Thus, $I = (R:J)
= (R:J^\star) = I_{\alpha_1} \cap I_{\alpha_2} \cap \ldots \cap
I_{\alpha_n}$, and so the only $M_\alpha$'s containing $I$ are
$M_{\alpha_1}, M_{\alpha_2}, \ldots, M_{\alpha_n}$. Since
$\bigcap_\alpha M_\alpha^2 \subseteq I$, we conclude that
$\{M_\alpha\} = \{M_{\alpha_1}, M_{\alpha_2}, \ldots,
M_{\alpha_n}\}$ is a finite family.

(ii)$\Rightarrow$(iii) First we show that $(J^\star:J^\star)R_M = (JR_M: JR_M)$,
for each nonzero ideal $J$ and for each $\star$-maximal ideal $M$ of $R$.
 For this, it is enough to show that $(JR_M: JR_M)  \subseteq
(J^\star:J^\star)R_M$. Let $x \in (JR_M: JR_M)$ and let $M_1, M_2,
\ldots, M_n$ be the $\star$-maximal ideals such that $xR_{M_i}
\neq R_{M_i}$. Since $R_{M_i}$ is
stable, for each $i=1,2,\ldots,n$, there exists  $y_i \in J$
such that $JR_{M_i} = y_i(JR_{M_i}:JR_{M_i})$ \cite[Lemma
3.1]{olb02}. Then, $xy_i
\in xJR_M \subseteq JR_M$ and there exists $d_i \in R
\smallsetminus M$ such that $d_ixy_i \in J$. Setting $d:= d_1
d_2 \ldots d_n $, we have $dxJR_{M_i} = dxy_i(JR_{M_i}:JR_{M_i})
\subseteq J(JR_{M_i}:JR_{M_i}) \subseteq JR_{M_i}$, for each $i =
1,2, \dots, n$. Moreover, if $N$ is a $\star$-maximal ideal such
that $N \not \in \{M_1,M_2, \ldots, M_n\}$, then $xR_N = R_N$.
Thus, $dxJR_N = dJR_N \subseteq JR_N$ for each $\star$-maximal
ideal $N$ of $R$ and so, $dxJ^\star = (dxJ)^\star = \bigcap (dxJR_N)
\subseteq \bigcap JR_N = J^\star$. It follows that $dx \in
(J^\star:J^\star)$ and $x \in (J^\star: J^\star)R_M$, since $d \in
R \smallsetminus M$. Thus, $(J^\star: J^\star)R_M = (JR_M:
JR_M)$.

Now let $T:=(J^\star: J^\star)$. We prove that there exists a finitely generated ideal
$H\sub J$ of  $R$ such that $(HT)^\star
= J^\star$. Let $N_1, N_2, \ldots N_s$ be the $\star$-maximal
ideals containing $J$. Since $R_{N_i}$ is stable, there exists $x_i \in J$, such that
$J^\star R_{N_i} = JR_{N_i} = x_i(J R_{N_i}: JR_{N_i}) = x_iTR_{N_i}$, for each $i =
1,2, \ldots, s$ \cite[Lemma
3.1]{olb02}. Let $F:= x_1R + x_2R + \ldots + x_sR$.
Since $FT\sub JT\sub J^\star$, we have
$JR_{N_i} = x_i T R_{N_i} \subseteq F T R_{N_i} \subseteq
J R_{N_i}$. It follows that $JR_{N_i} = FTR_{N_i}$ for
each $i = 1,2,\ldots, s$. If $F$ is not contained in any
$\star$-maximal ideal distinct from $N_1, N_2, \ldots, N_s$, we
have $JR_N = FTR_N$ for each $\star$-maximal ideal $N$ of
$R$ and so $J^\star = (FT)^\star$. Otherwise, let $N_{s+1},
N_{s+2}, \ldots, N_t$ be the $\star$-maximal ideals of $R$ containing
$F$ and not containing $J$. If $x \in J \smallsetminus (N_{s+1}
\cup N_{s+2} \cup \ldots \cup N_t)$ and  $H:= F+ xR$,
as before we get that $JR_N = HTR_N$ for each $\star$-maximal ideal $N$ of
$R$ and so $J^\star = (HT)^\star$.

(iii)$\Rightarrow$(i) We have to prove that $I:=
(J((J^\star:J):J))^\star = (J^\star:J^\star)$. By Lemma
\ref{lemma:2}, we have $I = \bigcap J((J^\star:J):J)R_M =
\bigcap JR_M ((J R_M:JR_M):JR_M) = \bigcap (JR_M:JR_M) =
\bigcap(J^\star:J^\star)R_M = (J^\star:J^\star)^\star = (J^\star:
J^\star)$ where the intersection varies over the set of all
$\star$-maximal ideals of $R$.

\smallskip
The fact that $\star = w$ is a straightforward consequence of
Corollary \ref{til=w}.
\end{proof}

We state explicitly the previous theorem for $\star = w$. A direct proof of (i)$\lra$(ii) is given in \cite[Theorem 2.2]{SS2}.

\begin{cor}\label{wstable} The following conditions are equivalent for an integral domain $R$:
\begin{enumerate}
\item[(i)] $R$ is $w$-stable.

\item[(ii)] $R$ has $t$-finite character and $R_M$ is stable,
for each $M\in t\op{-Max}(R)$.

\item[(iii)] $J^w$ is $\dot w$-finite in $(J^w:J^w)$, for each nonzero ideal $J$ of $R$, and $R_M$ is stable, for each $M\in t\op{-Max}(R)$.
 \end{enumerate}
\end{cor}

For $\star = d$, we obtain the following result,  
where (i)$\lra$(iii) is due to B. Olberding \cite[Theorem 3.3]{olb02}.

\begin{cor} \label{d=w}
The following conditions are equivalent for an integral domain $R$:
\begin{enumerate}
\item[(i)] $R$ is stable.

\item[(ii)] $R$ is $w$-stable and $d=w$.

\item[(iii)] $R$ has finite character and $R_M$ is stable,
for each maximal ideal $M$ of $R$.

\item[(iv)]  Each nonzero ideal $J$ of $R$ is finitely generated in $(J:J)$ and $R_M$ is stable
for each maximal ideal $M$ of $R$.
\end{enumerate}
In particular, a one-dimensional
$w$-stable domain is stable.
\end{cor}

We recall that a domain $R$ with the property that $\bigcap_{M\in\La_1} R_M
\neq \bigcap_{N\in\La_2} R_N$, for any two distinct subsets
$\La_1$ and $\La_2$ of $\op{Max}(R)$ is called a $\#$-domain. If $R$
has the same property for $\La_1,\, \La_2 \sub \tmax(R)$, we say that
$R$ is a \emph{$t\#$-domain} \cite{GHL}.

\begin{cor} \label{properties1} Let $R$ be a $w$-stable integral domain. Then:
\begin{enumerate}
\item[(1)] $\tspec(R)$ is treed.
\item[(2)] $R$ satisfies the ascending chain condition on $t$-prime ideals.
\item[(3)] $R$ is a $t\#$-domain.
\end{enumerate}
\end{cor}
\begin{proof} $(1)$ follows from Theorem \ref{thm:1} and the fact that
the spectrum of a local stable domain is linearly ordered \cite[Theorem
4.11(ii)]{olb02}.

$(2)$ follows from Theorem \ref{thm:1} and the fact that a local
stable domain satisfies the ascending chain condition on
 prime ideals \cite[Theorem 4.11(ii)]{olb02}.

$(3)$ follows from \cite[Corollary 1.3]{GHL}, because each $t$-maximal ideal of $R$ is divisorial by Lemma \ref{max}.
\end{proof}

\section{Overrings of $\star$-stable domains}

B. Olberding proved that overrings of stable domains are stable \cite[Theorem 5.1]{olb02}. This result was generalized by S. El Baghdadi, who showed that a $t$-linked overring $T$ of a $w$-stable domain is $w^\prime$-stable, where $w^\prime$ denotes the $w$-operation on $T$ (see \cite[Theorem 2.10]{SS2}).
Recall that an overring $T$ of an integral domain $R$ is called
\emph{$t$-linked} over $R$ if $T^ w =T$ \cite{DHLZ}. Each overring of $R$ is $t$-linked precisely when $d=w$ \cite[Theorem 2.6]{DHLZ}.
El Baghdadi's proof works more in general for semistar operations spectral and of finite type.

\begin{thm} \label{overring}
Let $R$ be a domain and $\star$ a semistar operation on $R$  such that $\star = \tilde{\star}$.  If $R$ is $\star$-stable, then each overring $T$ of $R$ is $\dot{\star}$-stable.
\end{thm}

\begin{proof}
Since $R\sub T$ implies $R^\star\sub T^\star$,
by Proposition \ref{rpunto} we can assume that $R^\star=R$ and $T^\star = T$, that
is, that $\star$ is a (semi)star operation on $R$ and $\dot{\star}$ is a (semi)star operation on $T$.

First we show that $T$ is $\dot{\star}$-locally stable. Let $M=M^{\dot{\star}}$ be
a $\dot{\star}$-maximal ideal of $T$. Then $(M \cap R)^\star \sub M^\star \cap R^\star = M \cap R$. Hence $M \cap R$ is a
$\star$-prime ideal of $R$ and $R_{M\cap R}\sub T_M$.
Since each localization of $R$ at
a $\star$-maximal ideal  is stable (Theorem \ref{thm:1}) and overrings of stable domains are stable \cite[Theorem 5.1]{olb02}, then  $T_M$ is stable.

In order to apply Theorem
\ref{thm:1}, we have to prove that $T$ has
$\dot{\star}$-finite character. Let $N$ be a $\star$-maximal ideal
of $R$ and let $\{M_\alpha\}$ be a family of
$\dot{\star}$-maximal ideals of $T$, such that $\bigcap_\alpha
M_\alpha \neq (0)$ and $M_\alpha \cap R \subseteq N$. We want to
show that $\{M_\alpha\}$ is a finite set. Let
$S:=\bigcap_\alpha T_{M_\alpha} \supseteq T$. Since $R_N \subseteq
T_{M_\alpha}$ for each $\alpha$, we have that $R_N \subseteq S$.
Hence $S$ is stable as an overring of the stable domain
$R_N$. Let $P_\alpha:= M_\alpha T_{M_\alpha} \cap S$, for each
$\alpha$. The $P_\alpha$'s are pairwise incomparable, because
$S_{P_\alpha} = T_{M_\alpha}$. Since $\bigcap_\alpha M_\alpha$ is
nonempty, also $\bigcap_\alpha P_\alpha$ is nonempty. Let $x \in
\bigcap_\alpha P_\alpha$. If the $P_\alpha$'s are infinitely many,
then $x$ is contained in infinitely many maximal ideals of $S$,
because $\spec(S)$ is treed \cite[Theorem 4.11(ii)]{olb02}. This
contradicts the finite character of $S$. It follows that the
$P_\alpha$'s, and so the $M_{\alpha}$'s, are finitely many. It is
easy to see that this implies the $\dot{\star}$-finite
character for $T$.
\end{proof}

\begin{cor} \label{dot}
Let $R$ be a $w$-stable domain and let $T$ be a $t$-linked overring of
$R$. Then $\dot{w} = w^\prime$ is the $w$-operation on $T$ and $T$ is $w'$-stable.
\end{cor}
\begin{proof}
$T$ is $\dot{w}$-stable by Theorem \ref{overring}. Since $\dot{w}$ is
a (semi)star operation on $T$ (because $T$ is $t$-linked over $R$)
and  $\dot{w} = \tilde{\dot{w}}$, it follows from
Corollary \ref{til=w} that $\dot{w} = \tilde{\dot{w}}=w^\prime$.
\end{proof}

If $\star$ is a semistar operation on $R$, the \emph{$\star$-integral
closure} of $R$ is the integrally closed overring of $R$ defined by
$R^{[\star]}:= \bigcup \{(F^\star:F^\star) \vert \, F\in \boldsymbol{f}(R) \}$ \cite{fl01}. Clearly $R^{[\star]}=R^{[\star_f]}$.
The $v$-integral closure $R^{[v]}$ of a domain $R$ is also called the \emph{pseudo-integral closure of $R$} \cite{AHZ}.
We say
that $R$ is \emph{$\star$-integrally closed} if $R^{[\star]}=R$.
In this case, it is easy to see that $\star$ is necessarily a
(semi)star operation on $R$ \cite[p. 50]{EFP}.  Denoting by $R^\prime$ the integral closure of $R$, we have $R\sub R^\prime\sub R^{[\star]}$. In addition, if $\star$ is a (semi)star operation and $\widetilde{R}:= \bigcup \{(I^v:I^v) \vert \, I\in \boldsymbol{F}(R) \}$ is the complete integral closure of $R$, we have $R^{[\star]}\sub R^{[v]}
\sub \widetilde{R}$.

It is known that the integral closure of a stable domain is a
Pr\"ufer domain \cite[Proposition 2.1]{rush}. We now show that, when $\star$ is
a semistar operation of finite type,
the $\star$-integral closure of a $\star$-stable domain is a
P$v$MD. Recall that an integrally closed domain is a P$v$MD if and only if $w=t$ \cite[Theorem 3.5]{K}.


\begin{thm} \label{prop:p*md}
Let $R$ be an integral domain and $\star$ a semistar operation on
$R$.  Assume that $R$ is $\star$-stable. Then, denoting by $t^\prime$ and $w^\prime$   respectively the $t$-operation and the $w$-operation on $R^{[\star]}$:

\begin{itemize}
\item[(1)] Each nonzero finitely
generated ideal of $R^{[\star]}$ is $\dot{\star}$-invertible.
\item[(2)]
If  $\star$ is of finite type,   then $R^{[\star]}$ is
a  P$v$MD and $\dot{\star} = t^\prime=w^\prime$.
\item[(3)]
If  $\star=\tilde{\star}$,   then $R^{[\star]}$ is
a  $w^\prime$-stable  P$v$MD and $\dot{\star} = t^\prime=w^\prime$.
\end{itemize}
\end{thm}
\begin{proof}

(1) Let $I$ be a nonzero finitely generated ideal of $R^{[\star]}$. We have to prove
that $(I(R^{[\star]}:I))^\star = R^{[\star]}$.

There exist $x \in
K$ and a finitely generated ideal $J$ of $R$, such that $I =
xJR^{[\star]}$. Since $R$ is $\star$-stable, $J^\star$ is
invertible in $(J^\star:J^\star) \subseteq R^{[\star]}$. Thus
$R^{[\star]} = (J^\star:J^\star)R^{[\star]} =
(J^\star((J^\star:J^\star):J))^\star R^{[\star]} \subseteq
(JR^{[\star]}(R^{[\star]}:JR^{[\star]}))^\star \subseteq
R^{[\star]}.$
It follows that $(I ( R^{[\star]}:I))^\star =
(JR^{[\star]}(R^{[\star]}:JR^{[\star]}))^\star =  R^{[\star]}$.

(2) If $\star$ is of finite type, $\dot{\star}$ is a (semi)star operation on
$R^{[\star]}$ \cite[Proposition 4.5(3)]{fl01}. By (1) and  \cite[Lemma
2.1]{chang/park}   we have
$I^{\dot{\star}} = I^{t^\prime}$ for each finitely generated ideal
$I$ of $R^{[\star]}$. Hence, $\dot{\star} = t^\prime$ and, again by (1),
we conclude that $R^{[\star]}$ is
a  P$v$MD.

(3) follows from (2) and Theorem \ref{overring}.
\end{proof}

\begin{cor} \label{staric}
Let $R$ be an integral domain and $\star$ a semistar operation of finite type
on $R$.  If $R$ is $\star$-stable and $\star$-integrally closed, then $R$ is a P$v$MD and $\star= t=w$.
In particular, a $t$-stable pseudo-integrally closed domain is a P$v$MD.
\end{cor}

\begin{cor} \label{krull}
Let $R$ be an integral domain and $\star$ a (semi)star operation of finite type
on $R$.  If $R$ is $\star$-stable and  completely integrally closed, then $R$ is a Krull domain and $\star= t=w$.
\end{cor}
\begin{proof} If $R$ is $\star$-stable, then each $t$-maximal ideal of $R$ is divisorial (Proposition \ref{max}). Hence $R$, being completely integrally closed, is a Krull domain  \cite[Theorem 2.6]{Gab}. Since $R$ is $\star$-integrally closed (because $R^{[\star]}
\sub \widetilde{R}$), then $\star= t=w$ (Corollary \ref{staric}).
\end{proof}

\begin{cor}\label{wic} \cite[Corollary 2.11]{SS2}
Let $R$ be a $w$-stable domain. Then the $w$-integral closure
$R^{[w]}$ of $R$ is a $w^\prime$-stable P$v$MD and $\dot{w}=t^\prime=w^\prime$  (where $t^\prime$ and $w^\prime$ are  respectively the $t$-operation and the $w$-operation on $R^{[w]}$) .
\end{cor}

\begin{rem} The integral closure $R^\prime$ of a $w$-stable domain is always $\dot{w}$-stable by Theorem \ref{overring}. However, $R^\prime$ is not necessarily $t$-linked over $R$ \cite[Section 4]{DHLRZ} and so  we cannot use Corollary \ref{dot} to conclude that $R^\prime$ is $w^\prime$-stable. (In fact,  when $\dot{w}$  is not a (semi)star operation, we cannot compare $\dot{w}$ and $w^\prime$.)
\end{rem}

 We do not know whether, for $\star$ of finite type, a $\star$-stable integrally closed domain is a P$v$MD. However,
when $\star$ is a
(semi)star operation on $R$ such that $\star = \tilde{\star}$, then $R$ is $\star$-integrally closed if and
only if $R$ is integrally closed \cite[Lemma 4.13]{EFP}. Hence, from Corollary  \ref{staric}, we have:

\begin{cor} Let $R$ be a domain and $\star$ a (semi)star operation on $R$  such that $\star = \tilde{\star}$.  If $R$ is $\star$-stable and integrally closed, then  $R$ is a P$v$MD and $\star= t=w$.
In particular, a $w$-stable integrally closed domain  is a P$v$MD.
\end{cor}

As a matter of fact, it is proved in \cite{SS2} that a
$w$-stable integrally closed domain is  precisely a strongly discrete P$v$MD with $t$-finite character.
Recall that a valuation domain $V$ is called \emph{strongly discrete} if
$PV_{P}$ is a principal ideal for each prime ideal $P$ of $V$. We
say that a P$v$MD (respectively, a Pr\"ufer domain) $R$ is {\it strongly discrete} if $R_{P}$ is a
strongly discrete valuation domain, for each $P\in
t\op{-Spec}(R)$ (respectively, for each $P\in \op{Spec}(R)$).
If $R$ is a strongly discrete P$v$MD (respectively, a Pr\"ufer domain)
and each proper $t$-ideal (respectively, each nonzero proper ideal) of
$R$ has
only finitely many minimal primes, then $R$ is called a
\emph{generalized Krull domain} \cite{elb1} (respectively a
\emph{generalized Dedekind domain}).

The following characterization of $w$-stable integrally closed domain is given in \cite[Theorem 2.6]{SS2}.
For stable domains, an analogous result is due to B. Olberding \cite{olb98, olb01}.

\begin{thm} \label{wstableic} The following conditions are
equivalent for an integral domain $R$:
\begin{enumerate}
    \item[(i)] $R$ is integrally closed and $w$-stable;
\item[(ii)] $R$ is a $w$-stable P$v$MD;
 \item[(iii)] $R$ is a strongly discrete P$v$MD with $t$-finite character;
 \item[(iv)] $R$ is a generalized Krull domain with $t$-finite character;
  \item[(v)] $R$ is a $w$-stable generalized Krull domain;
  \item[(vi)] $R$ is a P$v$MD with $t$-finite character and each
$t$-prime ideal of $R$ is $w$-stable;
 \item[(vii)] $R$ is $w$-stable and each $t$-maximal ideal of $R$ is
$t$-invertible.
\end{enumerate}\end{thm}

\section{$\star$-Divisorial $\star$-stable domains}

Following \cite{AZ1}, we say that a nonempty family $\Lambda$ of
nonzero prime ideals of $R$ is \emph{of finite character} if each
nonzero element of $R$ belongs to at most finitely many members of
$\Lambda$ and we say that $\Lambda$ is \emph{independent} if no
two members of $\Lambda$ contain a common nonzero prime ideal. We
note that, for a (semi)star operation $\star$ of finite type, the
family $\star\op{-Max}(R)$ is independent if and only if no two
members of $\star\op{-Max}(R)$ contain a common prime $t$-ideal,
because a minimal prime of a principal ideal is a $t$-ideal. When
the family of all maximal (resp. $t$-maximal) ideals of $R$ is
independent of finite character, $R$ is called an \emph{$h$-local
domain} (resp. a \emph{weakly Matlis domain}).

 A domain such that each ideal is divisorial (that is
$d=v$) is called a \emph{divisorial domain} and a domain whose
overrings are all divisorial is called \emph{totally divisorial}
\cite{BS}.
We say that a domain $R$ is \emph{$\star$-divisorial} if $\star = v$.

\begin{thm}\label{sstabledivisorial}
Let $R$ be an integral domain and $\star$ a (semi)star operation on
$R$. If $\star=\tilde{\star}$, the following conditions are equivalent:

\begin{enumerate}
\item[(i)] $R$ is $\star$-stable and $\star$-divisorial.

\item[(ii)] The family $\star\op{-Max}(R)$ is independent of finite character and
$R_M$ is totally divisorial, for each $M \in \star\op{-Max}(R)$.

\item[(iii)] Each overring $T$ of $R$ such that $T^\star = T$ is
$\dot{\star}$-divisorial.

\item[(iv)] $(I^\star:I^\star)$ is $\dot{\star}$-divisorial, for
each nonzero ideal $I$ of $R$.

\item[(v)] $v(I^\star) = \dot{\star}$ on $(I^\star:I^\star)$, for each nonzero ideal $I$ of $R$.
\end{enumerate}

Under these conditions, $\star = w$ on $R$ and
$\dot{\star}=\dot{w}=w^\prime$ is the $w$-operation on each $t$-linked
overring $T$ of $R$.
\end{thm}

\begin{proof}
(i)$\Rightarrow$(ii) Since $\star = \tilde{\star} \leq w \leq t\leq v$, if $R$ is
$\star$-divisorial, it is also $w$-divisorial. Hence the family
$\star\op{-Max}(R)=\tmax(R)$ is independent of finite character 
\cite[Theorem 1.5]{SS1}.

If $M$ is a $t$-maximal ideal of $R$,  $R_M$ is stable by Theorem
\ref{thm:1}. Hence, to show that  $R_M$ is totally divisorial,  it is enough to show that $R_M$ is divisorial  \cite[Theorem 3.12]{olb01}.
Let $I =JR_M$ be a nonzero ideal of $R_M$, with $J$ an ideal of $R$. Since $\tmax(R)$ is independent of finite character, we have
$I^v = (JR_M)^{v^\prime} = J^v R_M$, where $v^\prime$ is the $v$-operation on $R_M$ \cite[Corollary 5.3]{AZ1}.  Since
$J^\star = J^v$,  then $I^v = J^\star R_M = JR_M = I$. Hence,
$R_M$ is divisorial. 

(ii)$\Rightarrow$(iii) Let $T$ be an overring of $R$ such that
$T^\star=T$. By applying \cite[Corollary 5.2]{AZ1} for
$\mathcal{F} = \star\op{-Max}(R)$, if $J$ is a nonzero ideal of $T$ and
$M$ is a $\star$-maximal ideal, we have $(T:(T:J))R_M \subseteq
(TR_M: (T:J)R_M) = (TR_M:(TR_M:JR_M)) = JR_M$, where the last
equality holds because $R_M$ is totally divisorial and $TR_M
\supseteq R_M$. Thus, denoting by $v^\prime$ the $v$-operation of
$T$, we have $J^{v^\prime}R_M = JR_M$, for each $\star$-maximal
ideal $M$ of $R$. It follows that $(J^{v^\prime})^{\dot{\star}} =
J^{\dot{\star}}$. Since $\dot{\star}$ is a (semi)star operation on
$T$, \, $\dot{\star} \leq v^{\prime}$ and
$J^{v^\prime}=J^{\star^\prime}$.

(iii)$\Rightarrow$(iv) It is straightforward, since
$(I^\star:I^\star)^\star = (I^\star:I^\star)$.

(iv)$\Rightarrow$(v) Since $v(I^\star)$ is a (semi)star operation
on $(I^\star: I^\star)$, we have $v(I^\star) \leq v^\prime$, the
$v$-operation of $(I^\star:I^\star)$. Moreover, $I^\star$ is a
$\dot{\star}$-ideal, so by Lemma \ref{lemma:0}(2), $\dot{\star}
\leq v(I^\star)$. Thus, since by hypothesis $\dot{\star} =
v^\prime$, we have $\dot{\star} = v(I^\star)$.

(v)$\Rightarrow$(i) $R$ is clearly $\star$-divisorial, taking
$I=R$. That $R$ is $\star$-stable is a consequence of Lemma
\ref{lemma:1}(iii)$\Rightarrow$(i). In fact, if $I$ and $J$ are
ideals of $R$ such that $(I^\star:I^\star) = (J^\star: J^\star)$
it is clear that $v(I^\star)=v(J^\star)$, since
both coincide with $\dot{\star}$.

\smallskip
To finish, if $R$ is $\star$-stable, we have $\star = w$ by
Theorem \ref{thm:1} and $\dot{w}=w^\prime$ is the $w$-operation on each
$t$-linked overring $T$ of $R$ by Corollary \ref{dot}.
\end{proof}

We state explicitly the theorem for $\star = w$.

\begin{cor}\label{wstabledivisorial} The following conditions are equivalent for an integral domain $R$:
\begin{enumerate}

\item[(i)] $R$ is $w$-stable and $w$-divisorial.

\item[(ii)] $R$ is a weakly Matlis domain and $R_M$ is totally divisorial, for each $M\in \tmax(R)$.

\item[(iii)] Each t-linked overring $T$ of $R$ is
$\dot{w}$-divisorial.

\item[(iv)] $(I^w:I^w)$ is $\dot{w}$-divisorial, for each nonzero ideal $I$
of $R$.

\item[(v)] $v(I^w)=\dot{w}$ on $(I^w:I^w)$, for each nonzero ideal $I$ of $R$.
\end{enumerate}

Under these conditions, $\dot{w}=w^\prime$ is the $w$-operation on each $t$-linked overring $T$ of $R$.
\end{cor}

Since each $t$-linked overring of a $w$-stable domain is $w^\prime$-stable (Corollary \ref{dot}), we get:

\begin{cor} \label{overwstabledivisorial}
Each $t$-linked overring of a $w$-stable $w$-divisorial domain is a $w^\prime$-stable $w^\prime$-divisorial domain.
\end{cor}

We do not know if in general the
condition that each $t$-linked overring $T$ of $R$ is  $w^\prime$-divisorial implies that $R$ is $w$-stable. However, we now show that this is true in the integrally closed case and in the Mori case.

The following result follows from Corollary
\ref{wstabledivisorial} and the fact that a valuation domain $V$
is totally divisorial if and only if it is strongly discrete
\cite[Proposition 7.6]{BS}. (Compare it with Theorem
\ref{wstableic}.) The equivalences (ii)$\lra$(iii)$\lra$(iv) are
part of \cite[Theorem 3.5]{SS1}. (i)$\lra$(ii) is \cite[Corollary 2.8]{SS2}.

\begin{cor} \label{icwstabledivisorial} The following conditions are equivalent for an integral domain $R$:
\begin{enumerate}

\item[(i)] $R$ is an integrally closed $w$-stable $w$-divisorial domain.

\item[(ii)] $R$ is integrally closed and each $t$-linked overring of $R$
is ${w^\prime}$-divisorial.

\item[(iii)] $R$ is a strongly discrete P$v$MD and a weakly Matlis
domain.

\item[(iv)] $R$ is a $w$-divisorial generalized Krull domain.
\end{enumerate}
\end{cor}

\begin{cor} \label{cic} Let $R$ be an integral domain. The following are
equivalent:

\begin{enumerate}
\item[(i)] $R$ is completely integrally closed and $w$-stable.
\item[(ii)] $R$ is is completely integrally closed and $w$-divisorial.
\item[(iii)] $R$ is a Krull domain.
\end{enumerate}
\end{cor}

\begin{proof}(i)$\lra$(iii) follows from Corollary \ref{krull}.
(ii)$\lra $(iii) was proved in \cite[Proposition 3.7]{SS1}.
\end{proof}

Mori domains whose $t$-linked overrings are all $w^\prime$-divisorial were studied in \cite{SS1}.
A Mori domain is $w$-divisorial if and only if $R_M$ is a divisorial Noetherian domain,
necessarily one-dimensional, for each $t$-maximal ideal $M$ \cite[Theorem 4.5 and Proposition 4.1]{SS1}.

\begin{cor}\label{Mori2} The following conditions are equivalent for an integral domain $R$:
\begin{itemize}

\item[(i)] $R$ is a Mori $w$-stable $w$-divisorial domain.

\item[(ii)] $R$ is a Mori domain and each $t$-linked overring is
$w^\prime$-divisorial.

\item[(iii)] $R$ is a Mori domain and $R_M$ is totally divisorial,
for each  $M\in \tmax(R)$.

\item[(iv)] $R$ has $t$-dimension one and each $t$-linked
overring of $R$ is $w^\prime$-divisorial.

\item[(v)] $R$ has $t$-dimension one and is $w$-stable and $w$-divisorial.

\item[(vi)] For each nonzero ideal $I$ of $R$, there are $a$,
$b\in R$ such that $I^w= (aR+bR)^ w$.
\end{itemize}

\end{cor}
\begin{proof} (i)$\ra$(ii) and (v)$\ra$(iv) by Corollary \ref{wstabledivisorial}(i)$\ra$(iii).

(ii)$\lra$(iii)$\lra$(iv)$\lra$(vi) are proved in \cite[Theorem
4.11]{SS1}.

We show that (ii) and (iii) imply (i). In fact, by (ii) $R$ is
$w$-divisorial and so weakly Matlis, and by (iii) $R_M$ is totally
divisorial for each $t$-maximal $M$. So, we can conclude by
applying Corollary \ref{wstabledivisorial}(ii)$\ra$(i). In the
same way, we get that (iii) and (iv) imply (v).
\end{proof}

For $\star = d$ we recover the following characterizations of
totally divisorial domains: (i)$\lra$(iii)$\lra$(iv) are due to
B. Olberding \cite[Theorem 3.12 and Corollary 3.13]{olb01},
(iv)$\lra$(v)$\lra$(vi) are due to G. Picozza \cite[Theorem
2.57]{P}.

\begin{cor} \label{olb+gp}
The following conditions are equivalent for an integral domain $R$:
\begin{enumerate}

\item[(i)] $R$ is stable and divisorial.

\item[(ii)] $R$ is stable and $w$-divisorial.

\item[(iii)] $R$ is h-local and $R_M$ is totally divisorial, for each maximal ideal $M$ of $R$.

\item[(iv)] $R$ is totally divisorial.

\item[(v)] $(I:I)$ is divisorial for each nonzero ideal $I$ of $R$.

\item[(vi)] Each nonzero ideal $I$ of $R$ is $m$-canonical in $(I:I)$.

\end{enumerate}
\end{cor}

For ease of reference, we state Corollary \ref{olb+gp} in the
integrally closed case \cite{olb98, olb01} and in the Noetherian case
\cite{BS} (see also \cite[Corollary 3.6 and Proposition 4.8]{SS1}).

\begin{cor} \label{icstabledivisorial} The following conditions are equivalent for an integral domain $R$:
\begin{enumerate}

\item[(i)] $R$ is an integrally closed stable divisorial domain.

\item[(ii)] $R$ is an integrally closed totally divisorial domain.

\item[(iii)] $R$ is an $h$-local strongly discrete Pr\"ufer domain.
\item[(iv)] $R$ is a divisorial generalized Dedekind domain.
\end{enumerate}
\end{cor}

\begin{cor} \label{Noeth} The following conditions are equivalent for an integral domain $R$:
\begin{enumerate}

\item[(i)] $R$ is a Noetherian stable divisorial domain.

\item[(ii)] $R$ is a Noetherian totally divisorial domain.

\item[(iii)] $R$ is a one-dimensional totally divisorial domain.

\item[(iv)] $R$ is 2-generated (that is, each ideal of $R$ is generated by 2 elements).
\end{enumerate}
\end{cor}

\begin{rem}
(1) It is easy to check that the overrings of $R$ of the form $(I:I)$, where $I$ is a nonzero ideal of $R$, are precisely the overrings with nonzero conductor in $R$.
But, if $R$ is totally divisorial, it is not true in general that all
the overrings of $R$ have nonzero conductor (that is,
totally divisorial domains are not always \emph{conducive}). For example,
any conducive totally divisorial Pr\"ufer domain is a strongly discrete
valuation domain \cite[Theorem 4.7]{giampa},
while Corollary \ref{icstabledivisorial} shows that there exist plenty of non-quasilocal totally
divisorial Pr\"ufer domains.

(2) A Noetherian stable domain is always
one-dimensional \cite[Proposition 2.1]{sv}. We do not know if a Mori
$w$-stable domain need to have $t$-dimension one. However, we can say
that it has $t$-dimension at most equal to 2. In fact, by Corollary
\ref{wstable}, a Mori domain $R$ is $w$-stable if and only if $R_M$ is
stable for each $t$-maximal ideal $M$. In addition,
it is known that if $P$ is a nonzero nonmaximal prime ideal of a stable
domain $R$, then $R_P$ is a strongly discrete valuation domain
\cite[Theorem 4.11]{olb02}. Since the Mori property localizes and a
Mori valuation domain is a DVR, we see that a quasilocal stable Mori
domain has dimension at most equal to 2.

(3) Examples of Mori $w$-stable $w$-divisorial domains that are
neither Noetherian nor Krull can be constructed by means of pullbacks,
as in \cite[Example 4.13]{SS1}. We do not know any example of a
one-dimensional stable Mori domain that is not Noetherian.
\end{rem}

The next theorem shows that the study of $w$-stable
$w$-divisorial domains can be reduced to the case where the domain $R$
has $t$-dimension at least equal to two and has no $t$-invertible
$t$-prime ideals. If $\La$ is a set of prime ideals of $R$, we set $R_{\F(\La)}:=\bigcap_{P\in \La} R_P$.

\begin{thm}\label{decomposition} Assume that $R$ is a $w$-stable $w$-divisorial domain.
Let $\La_1$ be the set of the $t$-invertible $t$-maximal ideals of
$R$, $\La_2$ be the set of the height-one $t$-maximal ideals of
$R$ that are not $t$-invertible, $\La_3:=\tmax(R)\sm\{\La_1 \cup
\La_2\}$ and set $R_i:= R_{\F(\La_i)}$, for $i=1, 2, 3$. (If $\La_i= \emptyset$, set $R_i:= K$.) Then:
\begin{itemize}
\item[(1)] If $\La_1\neq \emptyset$, $R_1$ satisfies the equivalent
conditions of Corollary \ref{icwstabledivisorial}.

\item[(2)]If $\La_2\neq \emptyset$, $R_2$ satisfies the equivalent
conditions of Corollary \ref{Mori2} and has no $t^\prime$-invertible
$t^\prime$-prime ideals (where $t^\prime$ is the $t$-operation on $R_2$).

\item[(3)] If $\La_3\neq \emptyset$, $R_3$ is a $w$-stable $w$-divisorial
domain of $t^\prime$-dimension strictly greater than one with no $t^\prime$-invertible
$t^\prime$-prime ideals (where $t^\prime$ is the $t$-operation on $R_3$).

\item[(4)] $R=R_1\cap R_2\cap R_3$.
\end{itemize}
\end{thm}

\begin{proof} Let $\La$ be a nonempty set of $t$-maximal ideals of $R$ and $T:=R_{\F(\La)}$.
Since $\tmax(R)$ has $t$-finite character, then
$t^\prime\op{-Max}(T)=\{PR_P\cap T\,; P\in \La\}$ \cite[Proposition 1.17]{G}. In addition,
for $M = PR_P\cap T \in t^\prime\op{-Max}(T)$, we have $T_M=R_P$. Recalling that an ideal of a domain with $t$-finite character is $t$-invertible if and only if it is $t$-locally principal, we get that $M$ is $t^\prime$-invertible in $T$ if and only if $P$ is $t$-invertible in $R$.

Since $T$ is $t$-linked over $R$, $T$ is $w^\prime$-stable and $w^\prime$-divisorial (Corollary \ref{overwstabledivisorial}).
Hence  (3) and (4) follow easily.

(1) If $P\in \La_1$, then $PR_P$ is principal. Hence $R_P$ is a stable quasi-local domain (Theorem \ref{thm:1}) with principal maximal
ideal; whence it is a valuation domain \cite[Lemma 4.5]{olb02}. It
follows that $R_1$ is an integrally closed $w^\prime$-stable domain and then it satisfies the
equivalent conditions of
Corollary \ref{icwstabledivisorial}.

(2) $R_2$ has $t^\prime$-dimension one and is $w^\prime$-stable and $w^\prime$-divisorial. Hence
$R_2$ satisfies the equivalent conditions of Corollary
\ref{Mori2}. Since the $t$-maximal ideals in $\La_2$ are not
$t$-invertible, then $R_2$ has no $t^\prime$-invertible prime ideals.
\end{proof}

\begin{cor} Let $R$ be a totally divisorial domain. With the
notation of the previous theorem:
\begin{itemize}
\item[(1)] If $\La_1\neq \emptyset$,  $R_1$ satisfies the equivalent
conditions of Corollary \ref{icstabledivisorial}.
 \item[(2)] If $\La_2\neq \emptyset$, $R_2$ satisfies the equivalent
conditions of Corollary \ref{Noeth} and has no invertible prime
ideals.
\item[(3)] If $\La_3\neq \emptyset$,  $R_3$ is a totally divisorial domain of dimension strictly greater
than one with no invertible prime ideals.
\end{itemize}
\end{cor}

\section{$v$-Coherence}

A domain is \emph{coherent} if the intersection of
any two finitely generated ideals is finitely generated.
B. Olberding proved that
a stable divisorial domain is coherent
\cite[Lemma 3.2]{olb01}, even though there are stable domains that are
not coherent \cite[Section 5]{olb01-2}.

We next show that $w$-stable $w$-divisorial domains are
$v$-coherent. Recall that $R$ is \emph{$v$-coherent} if the
intersection of any two $v$-finite ideals is $v$-finite; this is
equivalent to say that the ideal $(R:I)$ is $v$-finite for each
nonzero finitely generated ideal $I$ of $R$.
The class of $v$-coherent domains properly includes P$v$MD's, Mori domains
and coherent domains \cite{GH1}.
A divisorial $v$-coherent domain is coherent.

The following lemma is probably known; for completeness we include the proof.

\begin{lemma}\label{lcoh} An integral domain $R$ with $t$-finite
character is $v$-coherent if and only if $R_M$ is $v$-coherent, for
each $t$-maximal ideal $M$ of $R$.
\end{lemma}

\begin{proof} If $R$ is a a $v$-coherent domain, then each
generalized ring of quotients of $R$ is $v$-coherent \cite[Proposition 3.1]{G}.

Conversely, let $J$ be a finitely generated nonzero ideal of $R$. If $J^ v
\neq R$, there are just finitely many $t$-maximal ideals $M_1, \dots,
M_n$ containing $J$ and, for each $i=1, \dots, n$, there is a
finitely generated ideal $H_i\sub (R:J)$ such that
$(R:J)R_{M_i} = (R_{M_i}:JR_{M_i}) = (H_iR_{M_i})^{v^\prime} = H_i^vR_{M_i}$ (where $v^\prime$ is the $v$-operation on $R_M$).

Let $H:= H_1+ \dots+ H_n$. If  $(R:J)\neq H^v$, let $N_1, \dots, N_m$ be the
$t$-maximal ideals of $R$ different from the $M_i$'s such that
$HR_{N_j} \neq R_{N_j}$, $j=1, \dots, m$. If $x \in R \sm \{N_1
\cup \dots \cup N_m\}$, by cheking $t$-locally, we get $(R:J)=
(H+xR)^v$.
 \end{proof}

\begin{thm} \label{vcoh} A $w$-stable $w$-divisorial domain is $v$-coherent.
\end{thm}
\begin{proof}
By Corollary \ref{wstabledivisorial} and Lemma \ref{lcoh}, because
totally divisorial domains are coherent \cite[Lemma 3.2]{olb01}.
\end{proof}

For $d=w$, the following proposition recovers \cite[Lemma 4.1]{olb01-2}.

\begin{prop} \label{finmax} A $w$-stable domain $R$ is $v$-coherent if and only if
each $t$-maximal ideal of $R$ is $v$-finite.
\end{prop}

\begin{proof}
Assume that $R$ is $v$-coherent and let $M\in \tmax(R)$. Since $M$ is
divisorial (Lemma \ref{max}), then $M= xR\cap R=(R: R+x^{-1}R)$, for some $x\in
R$. Thus $M$ is $v$-finite.

Conversely, if $R$ is $w$-stable, $R$ has $t$-finite character and
$R_M$ is stable, for each $t$-maximal ideal $M$ of $R$ (Theorem
\ref{thm:1}). Thus, $MR_M$ is divisorial in $R_M$. In addition, if
$M=J^v$, for some finitely generated ideal $J$ of $R$, then
$MR_M=J^vR_M=(JR_M)^{v^\prime}$ is $v^\prime$-finite (where $v^\prime$ is the $v$-operation on $R_M$). Hence, by Lemma \ref{lcoh},
we may assume that $R$ is stable quasilocal and that its maximal
ideal $M$ is divisorial $v$-finite.

We have to show that $(R:I)$ is $v$-finite for each nonzero finitely
generated ideal $I$ of $R$. If $I$ is ($t$-)invertible, this is true.
Thus
we may assume that $(R:I)=(M:I)$. Now, $(I:I)(M:I)=(M:I)=(M:M)(M:I)$
and, since $R$ is stable quasilocal, there exist $x, y \in K$ such that  
$I=x(I:I)$ and $M=y(M:M)$ \cite[Lemma 3.1]{olb02}. 
Hence, setting $\al :=xy$, $IM(M:I)=\al(M:I)=\al(R:I)$. On the other hand, we have also
$(R:I)=\be((R:I):(R:I))$, for some $\be \in K$. 
Thus we get
$\al^{-1}I^v=(R:IM(M:I))= (R:IM(R:I))
=(((R:I):(R:I)):M)=\be^{-1}(R:IM)$. It follows that
$(R:I)= \al^{-1}\be(IM)^v$ is $v$-finite.
\end{proof}

\begin{rem} (1) A one-dimensional stable coherent domain is
Noetherian, because its prime ideals are finitely generated. We do not
know whether a one-dimensional stable $v$-coherent domain must be Mori (or
even Noetherian).

(2) Generalized rings of quotients of $v$-coherent
domains are $v$-coherent \cite[Proposition 3.1]{G}, but it is not
known whether $t$-linked overrings of $v$-coherent domains are
$v$-coherent. However,
a $t$-linked overring of a $w$-stable $w$-divisorial
domain, being $w^\prime$-stable and $w^\prime$-divisorial (Corollary
\ref{overwstabledivisorial}), is $v$-coherent.
\end{rem}

By using $v$-coherence, we can see that $w$-stable
$w$-divisorial domains (respectively, totally divisorial domains) share
several properties with generalized Krull domains (respectively,
generalized Dedekind domains). In fact, since an integrally closed
$w$-stable $w$-divisorial domain (respectively, totally divisorial
domains) $R$ is a strongly discrete P$v$MD (respectively, Pr\"ufer
domain) (Theorem \ref{wstableic}), in the integrally closed case each one of
these properties becomes equivalent to $R$ being a generalized Krull
domain (respectively, generalized Dedekind domain) (see \cite[Theorems
3.5, 3.9 and Lemma 3.7]{elb1} and \cite[Corollary 2.15]{GHL}).

Recall that an overring $T$ of $R$ is said to be {\it $t$-flat} over
$R$ if $T_M=R_{M \cap R}$, for each $t$-maximal ideal $M$ of $T$
\cite{KP}. Flatness implies $t$-flatness, but the converse is not true
\cite[Remark 2.12]{KP}.  

\begin{cor}\label{properties2} Let $R$ be a $w$-stable $w$-divisorial domain. Then:
\begin{enumerate}
\item[(1)] Each $t$-prime ideal $P$ of $R$ is the radical of a $v$-finite divisorial ideal.
\item[(2)] $R$ satisfies the ascending chain condition on radical $t$-ideals.
\item[(3)] Each proper $t$-ideal has only finitely many minimal ($t$-)primes.
\item[(4)] For each $\La \sub \tspec(R)$, the overring $R_{\F(\La)}:=\bigcap_{P\in \La} R_P$ is a $t$-flat $t\#$-domain.
\end{enumerate}
\end{cor}

\begin{proof} (1) Since $R$ is $v$-coherent (Theorem \ref{vcoh}),
then $R_P$ is $v$-coherent, $w$-stable (Corollary \ref{dot}) and divisorial (Corollary \ref{wstabledivisorial}). Hence
$PR_P$ is $v^\prime$-finite in $R_P$(Proposition \ref{finmax}) and so $PR_P=(JR_P)^{v^\prime} = J^vR_P$, for some finitely
generated ideal $J$ of $R$. Since $\tmax(R)$ is independent of finite character (Corollary \ref{wstable}) and
$\tspec(R)$ is treed (Corollary \ref{properties1}), the set of minimal primes of $J^v$ is finite. Set
Min$(J^v)= \{P=P_1, \dots, P_n\}$. If $n\geq 2$, let $x\in P\sm (P_2\cup
\dots\cup P_n)$ and $I=(J+xR)^v$. Then $P=\sqrt I$.

(1) and (2) are equivalent by \cite[Lemma 3.7]{elb1}.

(3) follows from (1) by \cite[Lemma 3.8]{elb1}.

(4) Each localization of $R$ at a
$t$-prime ideal is divisorial  by Corollary \ref{wstabledivisorial} . In addition, $\tspec(R)$ is treed and
$R$ satisfies the ascending chain condition on $t$-prime ideals
by Corollary \ref{properties1}. Hence the overring $T := R_{\F(\La)}$ is $t$-flat by \cite[Corollary 2.12]{SS1}. Since $T$ is
$t$-linked over $R$, $T$ is $w^\prime$-stable and $w^\prime$-divisorial (Corollary
\ref{overwstabledivisorial}). Thus $T$ is a $t\#$-domain by Corollary
\ref{properties1}.
\end{proof}

\begin{cor} Let $R$ be a totally divisorial domain. Then:
\begin{enumerate}
\item[(1)] Each prime ideal of $R$ is the radical of a finitely
generated ideal.

\item[(2)] $R$ satisfies the ascending
chain condition on radical ideals.

\item[(3)] Each proper ideal
has only finitely many minimal primes.

\item[(4)] For each set $\La$ of prime ideals, the overring $R_{\F(\La)}:=\bigcap_{P\in \La} R_P$
 is a flat $\#$-domain.
\end{enumerate}
\end{cor}


\begin{thebibliography}{10}

\bibitem{AZ1}
D. D. Anderson and M. Zafrullah, \emph{Independent locally-finite
intersections of localizations}, Houston J. Math. \textbf{25}
(1999), no.~3, 433--452.

\bibitem{AHP}  D. D. Anderson, J. A. Huckaba and I. J. Papick, \emph{A note on stable domains}, Houston J. Math,

\bibitem{AHZ}  D. F. Anderson, E. Houston and M. Zafrullah. \emph{Pseudo-integrality}, Canad. Math. Bull. \textbf{34} (1991), 15--22.

\bibitem{b} H. Bass, \emph{On the ubiquity of Gorenstein rings}, Math Z. \textbf{82} (1963), 8--28

\bibitem{BS} S. Bazzoni and L. Salce, \emph{Warfield Domains}, J.
Algebra \textbf{185} (1996), 836--868.

\bibitem{chang/park}
G. W. Chang and J. Park, \emph{Star-invertible ideals of integral
domains}, Boll. Unione Mat. Ital. Sez. B Artic. Ric. Mat. (8)
\textbf{6} (2003), no.~1, 141--150.

\bibitem{chang/zafrullah} G.W. Chang and M. Zafrullah, \emph{The
$w$-integral closure of integral domains}, J. Algebra, to appear.

\bibitem{DHLZ} D. E. Dobbs, E. G. Houston, T. G. Lucas and M. Zafrullah,
\emph{$t$-linked overrings and Pr\"ufer $v$-multiplication
domains}, Comm. Algebra, \textbf{17} (1989), 2835--2852.

\bibitem{DHLRZ} D. E. Dobbs, E. G. Houston, T. G. Lucas, , M. Roitman and M. Zafrullah,
\emph{On $t$-linked overrings}, Comm. Algebra, \textbf{20} (1992), 1463--1488.

\bibitem{elb1} S. El Baghdadi, \emph{On a class of Pr\"ufer
$v$-multiplication domains}, Comm. Algebra \textbf{30} (2002),
3723--3742.

\bibitem{EFP} S.~El~Baghdadi, M. Fontana and G.~Picozza, \emph{Semistar
Dedekind domains}, J. Pure Appl. Algebra \textbf{193} (2004), no. 1-3,
27--60.

\bibitem{SS1}
S. El~Baghdadi and S. Gabelli, \emph{$w$-divisorial domains}, J. Algebra \textbf{285} (2005),
335--355.

\bibitem{SS2}S. El~Baghdadi and S. Gabelli, \emph{Ring-theoretic properties
of PVMDs}, submitted.

\bibitem{FR} D. Ferrand and M. Raynaud, \emph{Fibres formelles d'un
anneau local noeth\`erien}, Ann. Sci. \'Ecole Norm. Sup. \textbf{4}
(1970), 295--311

\bibitem{FH2000}
M. Fontana and J. Huckaba, \emph{Localizing systems and semistar
operations}, Non-Noetherian Commutative Ring Theory (Scott~T. Chapman and
Sarah Glaz, eds.), Kluwer Academic Publishers, 2000, pp.~169--198.

\bibitem{fjs}
M. Fontana, P. Jara and E. Santos, {\em Pr\"ufer
$\star$-multiplication domains and semistar operations}, J.
Algebra Appl. \textbf{2} (2003), 21--50.

\bibitem{fl01}
M. Fontana and K.~A. Loper, \emph{Kronecker function rings: a
general approach}, Ideal theoretic methods in commutative algebra
(Columbia, MO, 1999), Lecture Notes in Pure and Appl. Math., vol.
220, Dekker, New York, 2001, pp.~189--205.

\bibitem{fl03}
M. Fontana and K.~A. Loper, \emph{Nagata rings, {K}ronecker
function rings, and related semistar
operations}, Comm. Algebra \textbf{31} (2003), no.~10, 4775--4805.

\bibitem{Gab} S. Gabelli, \emph{Completely integrally closed domains and $t$-ideals},
Bollettino U. M. I. \textbf{7} (1989), 327--342.

\bibitem{G} S. Gabelli, \emph{On Nagata's Theorem for the class
group, II}, Lecture Notes in Pure and Appl. Math., vol. 206,
Marcel Dekker, New York, 1999, pp. 117--142.

\bibitem{GH1} S. Gabelli and E.G. Houston, \emph{Coherent-like conditions in
pullbacks}, Michigan Math. J. \textbf{44} (1997), 99--123.

\bibitem{GHL} S. Gabelli, E. Houston and T. Lucas, \emph{The
$t\#$-property for integral domains}, J. Pure Applied Algebra
\textbf{194} (2004), 281--298.

\bibitem{mcanI} W. J. Heinzer, J. A. Huckaba, and I.~J. Papick,
\emph{$m$-canonical ideals in integral domains}, Comm. Algebra
\textbf{26} (1998), no.~9, 3021--3043.

\bibitem{K} B. G. Kang, \emph{Pr\"ufer $v$-multiplication domains and the ring $R[X]_{N_v}$},
J. Algebra \textbf{123}, 151--170 (1989).

\bibitem{KP}D.J. Kwak and Y.S. Park, \emph{On $t$-flat overrings}, Chinese J.
Math. \textbf{23} (1995), 17--24.

\bibitem{L} J. Lipman, \emph{Stable ideals and Arf rings} J. Pure Applied Algebra \textbf{4} (1974),
319--336.

\bibitem{OM94}
A. Okabe and R. Matsuda, \emph{Semistar-operations on integral
domains}, Math. J. Toyama Univ. \textbf{17} (1994), 1--21.

\bibitem{olb98}
B. Olberding, \emph{Globalizing local properties of {P}r\"ufer
domains}, J. Algebra \textbf{205} (1998), no.~2, 480--504.

\bibitem{olb} B. Olberding, \emph{Stability of ideals and its applications}, Ideal theoretic methods in
commutative algebra
(Columbia, MO, 1999), Lecture Notes in Pure and Appl. Math., vol.
220, Dekker, New York, 2001, pp.~319--342.

\bibitem{olb01}
B. Olberding, \emph{Stability, duality, 2-generated ideals and a
canonical decomposition of modules}, Rend. Sem. Mat. Univ. Padova
\textbf{106} (2001), 261--290.

\bibitem{olb01-2}
B. Olberding, \emph{On the classification of stable domains},
J. Algebra \textbf{243} (2001), no.~1, 177--197.

\bibitem{olb02}
B. Olberding, \emph{On the structure of stable domains}, Comm.
Algebra \textbf{30} (2002), no.~2, 877--895.

\bibitem{P} G. Picozza, \emph{Semistar Operations and Multiplicative Ideal Theory},
Tesi di Dottorato, University  ``Roma Tre", Rome, June 2004.

\bibitem{giampa} G. Picozza
\emph {Star operations on overrings and semistar operations},
Comm. Algebra, \textbf{33} (2005), no.~6, 2051--203.

\bibitem{rush} D. E. Rush, \emph{Two-generated ideals and representations of
abelian groups over valuation rings}, J. Algebra \textbf{177}
(1995), 77--101.

\bibitem{sv} J. D. Sally and W. V. Vasconcelos, \emph{Stable rings}, J. Pure Appl. Algebra \textbf{4} (1974), 319--336.

\end{thebibliography}
\end{document}